\documentclass{article}
\usepackage{amssymb}
\usepackage{amsmath}

\begin{document}

\begin{center}
\bigskip 

{\Huge Variations on a Theme of Collatz}

\bigskip

{\Large Mario Bruschi and Francesco Calogero}

Physics Department, University of Rome "La Sapienza", Rome, Italy

Istituto di Fisica Nucleare, Sezione di Roma

\bigskip

\textit{Abstract}
\end{center}

Consider the recursive relation generating a new positive integer $n_{\ell
+1}$ from the positive integer $n_{\ell }$ according to the following simple
rules: if the integer $n_{\ell }$ is \textit{odd}, $n_{\ell +1}=3n_{\ell }+1$%
; if the integer $n_{\ell }$ is \textit{even}, $n_{\ell +1}=n_{\ell }/2$.
The so-called Collatz conjecture states that, starting from \textit{any
positive integer} $N$, the recursion characterized by the continued
application of these rules ends up in the cycle $4,$ $2,1$. This conjecture
is generally believed to be true (on the basis of extensive numerical
checks), but it is as yet unproven. In this paper---based on the assumption
that the Collatz conjecture is indeed true---we present a quite simple
extension of it, which entails the possibility to divide \textit{all natural
numbers} into \textit{three} disjoint classes, to each of which we
conjecture---on the basis of (not very extensive) numerical checks---that 
\textit{one third of all natural numbers} belong; or, somewhat equivalently,
to \textit{two} disjoint classes, to which we conjecture that respectively 
\textit{one third} and \textit{two thirds} of all natural numbers belong.

\bigskip

\section{Introduction}

The so-called "Collatz Conjecture" (CC) suggests that, starting from \textit{%
any positive integer} $N$, the "Collatz Recursion" (CR) stating that 
\begin{subequations}
\label{CR}
\begin{equation}
\text{if the integer }n_{\ell }\text{ is \textit{odd},~~~}n_{\ell
+1}=3n_{\ell }+1~,  \label{CRodd}
\end{equation}%
\begin{equation}
\text{if the integer }n_{\ell }\text{ is \textit{even},~~~}n_{\ell
+1}=n_{\ell }/2  \label{CReven}
\end{equation}%
ends up, after a \textit{finite} number of iterations, in an endless
repetition of the cycle $4,$ $2,1$.

This conjecture is generally believed to be true (on the basis of extensive
numerical checks: up to $N\approx 10^{20}$), but it is as yet unproven, and
the statement that "current mathematics is not up to provide a prove of it"
has been attributed to Paul Erd\"{o}s (26 March 1913--20 September 1996;
probably the most prolific and collaborative mathematician of all time, with
about 1500 papers published and about 500 co-authors). There is even some
justification to support the additional conjecture that the CC be an example
of \textit{unprovable mathematical truth}, the possible existence of such
propositions having been famously demonstrated by Kurt F. G\"{o}del (April
28, 1906-January 14, 1978); indeed in 1972 John. H. Conway has shown that a
natural generalization of the CC is \textit{algorithmically undecidable} 
\cite{C1972} \cite{KJ2007} \cite{L2010}. (For additional information on the
CC---including updates on its numerical checks---see the item "Collatz
Conjecture" in Wikipedia and the references indicated there).

In this paper---based on the assumption that the CC is valid---we discuss
some simple consequences of it, implying the possibility to \textit{separate
all natural numbers into three disjoint classes}, featuring a nontrivial
topology in the universe of all natural numbers. And we \textit{conjecture}%
---on the basis of (not very extensive) numerical checks---that these 
\textit{three} classes are of \textit{equal size} (with a precise definition
of this notion, see below). A simple modification of the CR is also
mentioned, and some of its implications.

\bigskip

\section{A variant of the Collatz Conjecture}

Consider the recursion relation that obtains by applying sequentially 
\textit{three times} the recursion CR, see (\ref{CR}); the resulting
recursion shall be hereafter indicated as CR3. It is then fairly obvious
that, if the CC is true, then, starting from \textit{any positive integer} $%
N $, repeated application of the recursion CR3 will eventually end in one of
its $3$ \textit{fixed points} $1,$ $2$ or $4;$ depending of course on the
value of the starting value $N$. This entails that the universe of natural
numbers $N$ gets divided in this manner into $3$ separate classes. A natural
question that comes to mind is the relative sizes of these $3$ classes. On
the basis of a rather amateurish numerical experimentation (see below), we
conjecture that \textit{one third of all the natural numbers} belong to each
one of these $3$ disjoint classes.

To make this notion precise, let us indicate with $\Sigma _{\lambda }\left(
S\right) $ the number of integers that, via the recursion RC3, end up,
respectively, in $\lambda =1,$ $2$, $4$ when starting from the integers $N$
in the interval $1\leq N\leq S$: implying of course 
\end{subequations}
\begin{equation}
\Sigma _{1}\left( S\right) +\Sigma _{2}\left( S\right) +\Sigma _{4}\left(
S\right) =S~.
\end{equation}%
A very simple experimentation---performed on a cheap PC---yielded the
following results: 
\begin{subequations}
\begin{equation}
\Sigma _{1}\left( 10^{5}\right) =33364~~~,\Sigma _{2}\left( 10^{5}\right)
=33311~~~,\Sigma _{4}\left( 10^{5}\right) =33325~,  \label{100K}
\end{equation}%
\begin{equation}
\Sigma _{1}\left( 10^{6}\right) =332858~~~,\Sigma _{2}\left( 10^{6}\right)
=333314~~~,\Sigma _{4}\left( 10^{6}\right) =333828~,  \label{1M}
\end{equation}%
\begin{equation}
\Sigma _{1}\left( 10^{7}\right) =3325705~~~,\Sigma _{2}\left( 10^{7}\right)
=3338680~~~,\Sigma _{4}\left( 10^{7}\right) =3335615~.  \label{10M}
\end{equation}%
This motivates us to proffer the following \textbf{Conjecture}: 
\end{subequations}
\begin{equation}
\underset{S\rightarrow \infty }{\lim }\left[ \frac{\Sigma _{\lambda }\left(
S\right) }{S}\right] =\frac{1}{3}~,~~~\lambda =1,~2,~4~.  \label{Conj}
\end{equation}

Actually, the idea to proffer this conjecture was based on the numerical
finding (\ref{100K}), which required only a few-minute computation; the
subsequent findings (\ref{1M}) and (\ref{10M}) damped our enthusiasm, but
nevertheless did not dissuade us from proffering the conjecture (\ref{Conj}%
), as a basis for further numerical checks (presumably to be performed by
less amateurish means).

\bigskip

\section{The \textit{paso doble} variant and its variant}

The CR implies that, if $n_{\ell }$ is \textit{odd}, then certainly $n_{\ell
+1}$ is \textit{even} (see (\ref{CR})). It is therefore quite natural to
introduce a variant of it which reads as follows: 
\begin{subequations}
\label{pdCR}
\begin{equation}
\text{if the integer }n_{\ell }\text{ is \textit{odd},~~~}n_{\ell +1}=\left(
3n_{\ell }+1\right) /2~,  \label{pdCRodd}
\end{equation}%
\begin{equation}
\text{if the integer }n_{\ell }\text{ is \textit{even},~~~}n_{\ell
+1}=n_{\ell }/2\text{~.}  \label{pdCReven}
\end{equation}%
This modified CR goes under various names in the literature; we take the
liberty to refer to it---hopefully for obvious reasons---as the "\textit{%
paso doble} Collatz Recursion" (\textit{pd}CR). Note that it does \textit{not%
} correspond to a double iteration of the CR recursion; indeed only (\ref%
{pdCRodd}) is different from (\ref{CRodd}), while (\ref{pdCReven}) is
identical to (\ref{CReven}).

It is now obvious that the CC implies that, starting from \textit{any
natural number }$N$, by iterating this new recursion \textit{pd}CR, see (\ref%
{pdCRodd}), one eventually ends up into the simplest cycle $1,2$; and again
this naturally suggest to introduce a new iteration---let us call it \textit{%
pd}CR2---each step of which amounts to \textit{two} steps of the iteration 
\textit{pd}CR. Clearly this iteration---starting from any arbitrary integer $%
N$---shall eventually settle on one of its \textit{two fixed points}, $1$ or 
$2;$ with twice as many integers $N$ eventually ending up on $1$ than on $2$
(if the conjecture (\ref{Conj}) is true). Implying another separation of 
\textit{all the natural numbers into two disjoint classes}; of course, a
separation which is quite trivially related to the subdivision in \textit{%
three} classes mentioned in the previous section.

\bigskip

\section{Envoy}

The above considerations are all based on quite elementary
mathematics---albeit relying on the validity of a conjecture the status of
which is still an \textit{open} mathematical issue. Of course further
investigations about the topology of the classes mentioned above are quite
open problems, as well as the amusing possibility to extend the CR from the 
\textit{natural numbers} to \textit{real} and even to \textit{complex numbers%
}; let us indeed confess that the discovery of this possibility was the main
motivation for our recent involvement in this kind of investigations
(although one of us was less new to this game than the other: see \cite%
{B20058}). But then we discovered that our discovery had been invented long
ago, as indeed reported in Wikipedia...

\end{subequations}

\end{document}